\documentclass[a4paper]{article}

\usepackage{RRA4}
\usepackage{hyperref}
\RRdate{Novembre 2006}
\RRauthor{% les auteurs
Monique Chyba
\thanks{University of Hawaii, Department of Mathematics, 2565 Mc Carthy the Mall, Honolulu 96822, Hawaii}%
  \and
Ernst Hairer
\thanks[gnv]{Universit\'e de Gen\`eve, Section de math\'ematiques, 2-4 rue du Li\`evre, CP 64, 1211 Gen\`eve 4, Switzerland}%
  \and
Gilles Vilmart
\thanksref{gnv}
\thanks{INRIA Rennes, projet IPSO, Campus de Beaulieu, 35042 Rennes Cedex, France}%
}
\authorhead{Chyba \& Hairer \& Vilmart}
\RRtitle{Intégrateurs symplectiques
en géométrie sous-riemannienne:\\ le problème de Martinet}
\RRetitle{Symplectic integrators
in sub-Riemannian geometry:\\ the Martinet case}
\titlehead{Symplectic integrators in sub-Riemannian geometry:
the Martinet case}
\RRnote{This work was partially supported by the 
Fonds National Suisse, project No.\ 200020-109158
and by the National Science Foundation grant DMS-030641.}
\RRresume{
On compare les performances d'intégrateurs symplectiques et non symplectiques
pour le calcul de géodésiques normales et de points conjugués dans un exemple 
sous-riemannien, le problème de Martinet. On étudie le problème d'abord avec une métrique
plate, puis avec une perturbation à un paramètre conduisant à des géodésiques
non intégrables.
De cette étude, on déduit que proche des directions anormales, une méthode
symplectique est bien plus efficace pour ce problème de contrôle optimal. 
L'explication repose sur la théorie de l'analyse rétrograde en intégration numérique
géométrique. 
}
\RRabstract{
We compare the performances of symplectic and non-symplectic
integrators for the computation of normal geodesics and conjugate
points in sub-Riemannian geometry at the example of the Martinet case.
For this case study we consider first the flat metric, and then a one 
parameter perturbation leading to non integrable geodesics. 
From our computations we deduce that near the abnormal directions
a symplectic method is much more efficient for this optimal control problem. 
The explanation relies on the theory of backward error analysis
in geometric numerical integration.
}
\RRmotcle{g\'eom\'etrie sous-riemannienne, Martinet, geod\'esique anormale, 
int\'egrateur symplectique, analyse r\'etrograde}
\RRkeyword{sub-Riemannian geometry, Martinet, abnormal geodesic,
symplectic integrator, backward error analysis.}
\RRprojet{Ipso}  % cas d'un seul projet
\RRtheme{\THNum} % cas d'un seul theme
 \URRennes  % pour ceux qui sont \`a l'ouest
\usepackage{amsmath}
\usepackage{epsfig}
\usepackage{GGGraphics}

\def\ts{\thinspace}
\def\rk2{\textsc{rk}$2$}
\def\verlet{St\"ormer--Verlet }
\def\bigo{\mathcal O}

\newtheorem{theo}{Theorem}[section]

\newtheorem{rmrk}[theo]{Remark}

\begin{document}
\makeRR   % cas d'un rapport de recherche
%% \makeRT % cas d'un rapport technique.
%% a partir d'ici, chacun fait comme il le souhaite
%%

\section*{Introduction}
This paper is a follow-up to a series of articles that were published 
in the past decades, see \cite{SR97,SR99} and the
references therein. 
There, the authors provide an analytic study of the singularity 
of the sub-Riemannian sphere in the Martinet case. They complement their 
analysis with numerical computations to represent the geodesics
and the sphere, and to locate the conjugate points. The integrator used for these 
computations is an explicit Runge-Kutta method of order 5(4). 
The goal of the present paper is to compare the performances of a symplectic integrator 
versus a non-symplectic one for this optimal control problem.  To be more precise, 
let $(U,\Delta,g)$ be a sub-Riemannian structure where $U$ is an open neighborhood 
of $R^3$, $\Delta$ a distribution of constant rank 2 and $g$ a Riemannian metric. 
When $\Delta$ is a contact distribution, there are no abnormal 
geodesics, and a non-symplectic integrator is as efficient as a symplectic one. 
However, when the distribution is taken as the kernel of the Martinet one-form, 
we show that a symplectic integrator is much more efficient
for the computation of the normal 
geodesics and their conjugate points near the abnormal directions. 

Both problems, the Martinet flat case and a non integrable 
perturbation, are introduced
in Sect.\ts \ref{sect:martinet} together with the corresponding
differential equations. Numerical experiments with
an explicit Runge--Kutta method and with the symplectic
\verlet method are presented in
Sect.\ts \ref{sect:compar} and illustrated with
figures. Close to abnormal geodesics,
the results are quite spectacular. 
For a relatively large step size, the symplectic integrator provides a solution with the correct qualitative
behavior and a satisfactory accuracy, while for the same step size the non-symplectic 
integrator gives a completely wrong numerical solution with an incorrect behavior,
particularly for the non integrable case. 
The explanation relies on the theory of backward error analysis (Sect.\ts \ref{sect:backward}).
It is related to the geometric structure of the problem and its solutions.

\section{A Martinet type sub-Riemannian structure} \label{sect:martinet}

In this section, we briefly recall some results of \cite{SR97}
for a sub-Riemannian structure
$(U,\Delta,g)$.
Here, $U$ is an open neighborhood of the origin
in $R^3$ with coordinates $q=(x,y,z)$,
and $g$ is a Riemannian metric
for which a graduated normal form, at order $0$, is
$g=(1+\alpha y)dx^2+(1+\beta x+\gamma y)dy^2$. The distribution $\Delta$
is generated by the two vector fields
$F_1=\frac{\partial}{\partial x}+\frac{y^2}{2}\frac{\partial}{\partial z}$ 
and $F_2=\frac{\partial}{\partial y}$
which correspond to $\Delta=\ker\omega$ where $\omega=dz-\frac{y^2}{2}dx$ is the Martinet 
canonical one-form. To this distribution
we associate the affine control system
\begin{equation*}
\dot q = u_1(t)F_1(q) +  u_2(t)F_2(q) ,
\end{equation*}
where $u_1(t), u_2(t)$ are measurable bounded functions which
act as controls. 

We consider two cases, the Martinet flat case $g=dx^2+dy^2$, 
an integrable situation, and a one parameter perturbation $g=dx^2+(1+\beta x)^2dy^2$ 
for which the set of geodesics is non integrable.

\subsection{Geodesics}

It follows
from the Pontryagin 
maximum principle, see \cite{SR97,SR99}, that the normal geodesics 
corresponding to $g=dx^2+(1+\beta x)^2dy^2$ are solutions of an Hamiltonian 
system
\begin{equation}
\label{eq:Hamiltonian}
\dot q=\frac{\partial H}{\partial p}(q,p), \qquad 
\dot p=-\frac{\partial H}{\partial q}(q,p) ,
\end{equation}
where $q=(x,y,z)$ is the state, $p=(p_x,p_y,p_z)$ is the adjoint state, 
and the Hamiltonian is
\begin{equation*}
H(q,p)=\frac{1}{2}\biggl(\Bigl(p_x+p_z\frac{y^2}{2}\Bigr)^2+\frac{p_y^2}{(1+\beta x)^2}\biggr).
\end{equation*}
In other words, the normal geodesics are solutions of the following equations:
\begin{equation}\label{eq:martinetsyst}
\begin{array}{rcl}
\dot x&=& \displaystyle p_x+p_z\frac{y^2}{2}\\[3mm]
\dot y&=& \displaystyle \frac{p_y}{(1+\beta x)^2}\\[3mm]
\dot z&=& \displaystyle \Bigl(p_x+p_z\frac{y^2}{2}\Bigr)\frac{y^2}{2}
\end{array}
\qquad
\begin{array}{rcl}
\dot p_x&=& \displaystyle \frac{\beta\, p_y^2}{(1+\beta x)^3}\\[3mm]
\dot p_y&=& \displaystyle -\Bigl(p_x+p_z\frac{y^2}{2}\Bigr)p_zy\\[3mm]
\dot p_z&=& \displaystyle 0  . \phantom{\frac{y^2}{2}}
\end{array}
\end{equation}

Notice that the variables $z$ and
$p_z$ do not influence the other equations
(except via the initial value $p_z(0)$), so that we are
actually confronted
with a Hamiltonian system in dimension four.
For the Martinet flat case ($\beta =0$), the interesting
dynamics takes place in the two-dimensional space of coordinates
$(y,p_y)$. The Hamiltonian is
\begin{equation*}
H(y,p_y)=\frac{p_y^2}{2}+ \frac 12 \Bigl(p_x+p_z\frac{y^2}{2}\Bigr)^2 ,
\end{equation*}
where $p_x$ and $p_z$ have to be considered as constants. This
is a one-degree of freedom mechanical system with a quartic
potential. For $p_x<0<p_z$, the Hamiltonian $H(y,p_y)$ has two local
minima at $(y\! =\! \pm \sqrt{-2p_x/p_z} ,\, p_y\! =\! 0)$, which correspond to
stationary points of the vector field. In this case, the origin
$(y\! =\! 0,\, p_y\! =\! 0)$ is a saddle point.

%\subsection{Abnormal geodesics}

Whereas normal 
geodesics correspond to oscillating motion,
it is shown in \cite{SR97,SR99} 
that the abnormal geodesics are 
the lines $z=z_0$ contained in the plane $y=0$.
For the considered metrics, the abnormal geodesics
can be obtained as projections of normal 
geodesics, we say that they are not strictly abnormal. 
In \cite{SR99}, the authors introduce a geometric framework to analyze the 
singularities of the sphere in the abnormal direction when $\beta \ne 0$. 
See also \cite{BLT2001,BT2001} for a precise description of the role of the 
abnormal geodesics in sub-Riemannian geometry in the general non-integrable case, 
i.e., when the abnormal geodesics can be strict. 
The major result of these papers is the proof that the sub-Riemannian sphere 
is not sub-analytic because of the abnormal geodesics.

Interesting phenomena arise when the normal geodesics
are close to the 
separatrices connecting the saddle point.
Therefore, we shall consider in Sect.\ts\ref{sect:compar}
the computation of normal geodesics with 
$y(0)=0$ and $p_y(0)>0$ but small.

\subsection{Conjugate points}

For the Hamiltonian system (\ref{eq:Hamiltonian})
we consider the exponential mapping
\begin{equation*}
\exp_{q_0,t}:p_0\longmapsto q(t,q_0,p_0)
\end{equation*}   
which, for fixed $q_0\in R^3$,
is the projection $q(t,q_0,p_0)$ onto the state space
of the solution of 
$(\ref{eq:Hamiltonian})$ starting at $t=0$ from $(q_0,p_0)$. 
Following the definition in \cite{SR97} we say that the point $q_1$ is conjugate
to $q_0$ along $q(t)$ if there exists $(p_0,t_1)$, $t_1>0$,
such that $q(t)=\exp_{q_0,t}(p_0)$ with
$q_1=\exp_{q_0,t_1}(p_0)$, and the mapping
$\exp_{q_0,t_1}$ is not an immersion at $p_0$. 
We say that $q_1$ is the first conjugate point if $t_1$ is minimal. 
First conjugate points play a major role when studying optimal control problems 
since it is a well known result that a geodesic is not optimal
beyond the first conjugate point. 

For the numerical computation of the first conjugate point,
we compute 
the solution of the Hamiltonian system  (\ref{eq:Hamiltonian}) together with 
its variational equation,
\begin{equation}\label{eq:variational}
\dot y = J^{-1}\nabla H(y), \qquad \dot \Psi = J^{-1}\nabla^2 H(y) \Psi .
\end{equation}
Here, $y=(q,p)$ and $J$ is the canonical matrix for Hamiltonian systems. 
It can be shown that for Runge-Kutta methods, the derivative of the numerical 
solution with respect to the initial value, $\Psi_n=\partial y_n/\partial y_0$, 
is the result of the same numerical integrator applied to the augmented system 
(\ref{eq:variational}), see \cite[Lemma VI.4.1]{hairer06gni}. 
Here, the matrix
\begin{equation*}
\Psi = \left(
\begin{array}{cc}
{\partial q}/{\partial q_0} & {\partial q}/{\partial p_0} \\
{\partial p}/{\partial q_0} & {\partial p}/{\partial p_0} 
\end{array}
\right)
\end{equation*}
has dimension $6\times6$.
The conjugate points are obtained when ${\partial q}/{\partial p_0}$ becomes 
singular, i.e., $\det({\partial q}/{\partial p_0})=0$.

\section{Comparison of symplectic and non-symplectic integrators}
\label{sect:compar}

For the numerical integration of the Hamiltonian system (\ref{eq:Hamiltonian}), 
where we rewrite $\frac{\partial H}{\partial q}(q,p)=H_q(q,p)$ and 
$\frac{\partial H}{\partial p}(q,p)=H_p(q,p)$, we consider two integrators of the same order $2$:
\begin{itemize}
\item[$\bullet$] a non-symplectic, explicit Runge--Kutta
discretization, 
denoted \rk2 (see \cite[Sect.\ts II.1.1]{hairer06gni}),
\begin{equation} \label{eq:rk2}
\begin{array}{rcl}
q_{n+1/2} &=& \displaystyle q_n + \frac h2 H_p(q_n,p_n) \\[3mm]
q_{n+1} &=& \displaystyle q_n + h H_p(q_{n+1/2},p_{n+1/2})
\end{array}
\qquad
\begin{array}{rcl}
p_{n+1/2} &=& \displaystyle p_n - \frac h2 H_q(q_n,p_n) \\[3mm]
p_{n+1} &=& \displaystyle p_n - h H_q(q_{n+1/2},p_{n+1/2})
\end{array}
\end{equation}
\item[$\bullet$] the symplectic \verlet scheme 
(see e.g. \cite[Sect.\ts VI.3]{hairer06gni}),
\begin{eqnarray}
p_{n+1/2} &=& p_n - \frac h2 H_q(q_n,p_{n+1/2}) \nonumber \\
q_{n+1} &=& q_n + \frac h2 \Bigl( H_p(q_n,p_{n+1/2}) 
+ H_p(q_{n+1},p_{n+1/2}) \Bigr) \label{eq:verlet}\\
p_{n+1} &=& p_{n+1/2} - \frac h2 H_q(q_{n+1},p_{n+1/2}) \nonumber 
\end{eqnarray}
\end{itemize}
where $q_n=(x_n,y_n,z_n)$ and $p_n=(p_{x,n},p_{y,n},p_{z,n})$. 
Here, $q_n \approx q(nh),\, p_n \approx p(nh)$ and $h$ is the step size. 

For the computation
of the conjugate points, we apply the numerical methods to
the variational equation~(\ref{eq:variational}).
Notice that only the partial derivatives with respect
to $p_0$ have to be computed.
Conjugate points are then
detected when $\det({\partial q_n}/{\partial p_0})$ 
changes sign. We approximate them by linear interpolation which
introduces an error of size
$\bigo(h^2)$. This is comparable to the accuracy of 
the chosen integrators which are both of second order.

\begin{rmrk} \rm
The \verlet scheme (\ref{eq:verlet})
is implicit in general. A few fixed point iterations yield the
numerical solution with the desired accuracy.
Notice however that the method becomes 
explicit in the Martinet flat case $\beta=0$. One simply has to compute 
the components in a suitable order, for instance 
$p_{x,n+1}$, $p_{z,n+1}$, $p_{y,n+1/2}$, $ y_{n+1}$, $x_{n+1}$, $
 z_{n+1}$, $p_{y,n+1}$. 
\end{rmrk}

\begin{figure}[b]
\centering
%\GGGinput{../prog/figs/}{figmartinet1}
 \begin{picture}(0,0)
  \epsfig{file=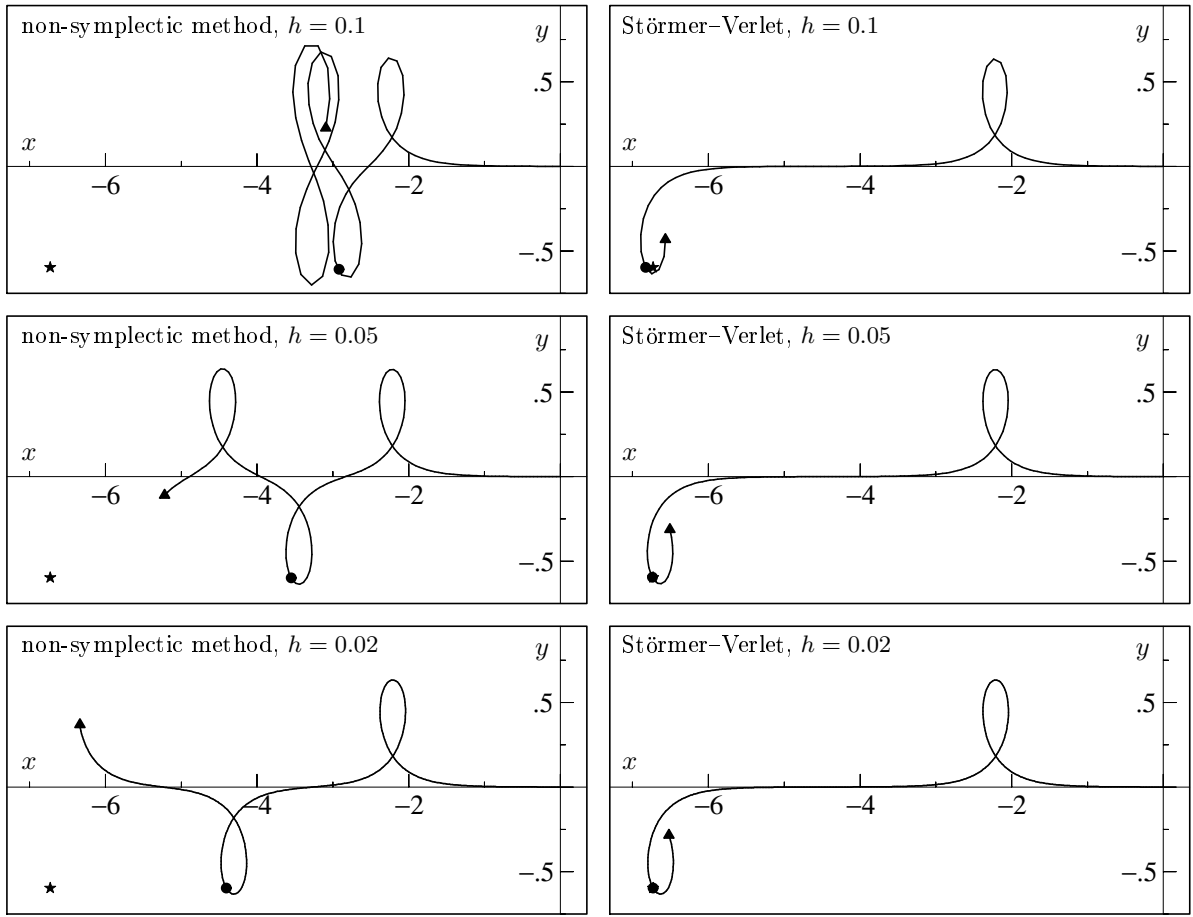}
 \end{picture}%
\begin{picture}(441.0,341.5)(  4.2,251.7)
  \GGGputW[   0, 100](  9.88,543.95){$x$}
  \GGGputW[   0, 100](202.12,586.83){$y$}
  \GGGputW[   0, 100](  9.88,590.08){\small non-symplectic method,  $h= 0.1$}
  \GGGputW[   0, 100](234.64,543.95){$x$}
  \GGGputW[   0, 100](426.87,586.83){$y$}
  \GGGputW[   0, 100](234.64,590.08){\small St\"ormer--Verlet,  $h= 0.1$}
  \GGGputW[   0, 100](  9.88,427.24){$x$}
  \GGGputW[   0, 100](202.12,470.24){$y$}
  \GGGputW[   0, 100](  9.88,473.37){\small non-symplectic method,  $h= 0.05$}
  \GGGputW[   0, 100](234.64,427.24){$x$}
  \GGGputW[   0, 100](426.87,470.24){$y$}
  \GGGputW[   0, 100](234.64,473.37){\small St\"ormer--Verlet,  $h= 0.05$}
  \GGGputW[   0, 100](  9.88,310.64){$x$}
  \GGGputW[   0, 100](202.12,353.52){$y$}
  \GGGputW[   0, 100](  9.88,356.77){\small non-symplectic method,  $h= 0.02$}
  \GGGputW[   0, 100](234.64,310.64){$x$}
  \GGGputW[   0, 100](426.87,353.52){$y$}
  \GGGputW[   0, 100](234.64,356.77){\small St\"ormer--Verlet,  $h= 0.02$}
 \end{picture}
\caption{Trajectories in the $(x,y)$-plane for the flat case $\beta=0$. 
\label{fig:figpaper1}}
\end{figure}

\subsection{Martinet flat case} \label{sect:flatcase}

We consider first the flat case $\beta =0$ in the
Hamiltonian system (\ref{eq:martinetsyst}).
As initial values we choose (cf.\ \cite{SR97})
\begin{equation}\label{eq:init}
x(0)=y(0)=z(0)=0, \quad 
p_x(0)=\cos \theta_0, \quad 
p_y(0)=\sin \theta_0, \quad p_z(0)=10, \quad 
\textrm{where} \quad 
\theta_0=\pi-10^{-3} ,
\end{equation}
so that we start close to an abnormal geodesics,
and we integrate the system over the interval
$[0,9]$.

Figure \ref{fig:figpaper1} displays the projection onto
the $(x,y)$-plane of the numerical solution obtained with different
step sizes~$h$ by the two integrators.
The initial value is at the origin, and the final state is
indicated by a triangle.
The circles represent the first conjugate point detected
along the numerical solution, while 
the stars give the position of the first conjugate point
on the exact solution of the problem. There is an enormous
difference between the two numerical integrators.
The symplectic (St\"ormer--Verlet) method (\ref{eq:verlet})
provides  a qualitatively correct solution 
already with a large step size $h=0.1$, and it gives an
excellent approximation for step sizes smaller than
$h=0.05$. On the other hand, the non-symplectic,
explicit Runge--Kutta method (\ref{eq:rk2}) gives completely
wrong results, and step sizes smaller than $10^{-3}$ 
are needed to provide an acceptable solution. An explanation
of the different behavior of the two integrators will be given
in Sect.\ts\ref{sect:backward} below.

\begin{figure}[ht]
\centering
%\GGGinput{}{figmartinet2}
%\GGGinput{../prog/figs/}{figmartinet2}
 \begin{picture}(0,0)
 \epsfig{file=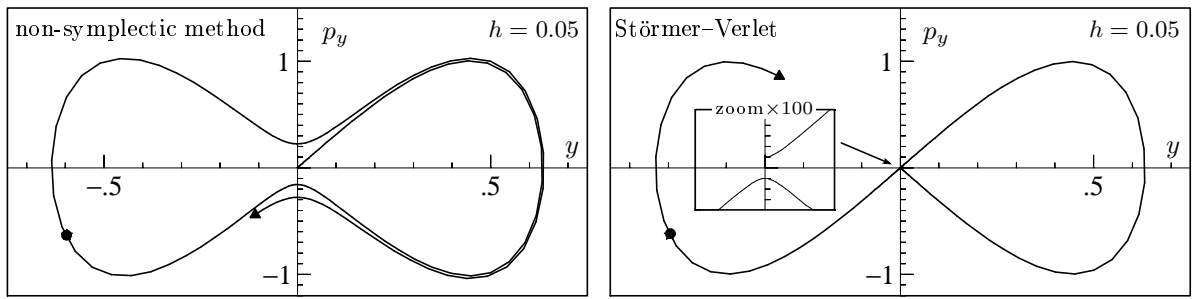}
 \end{picture}%
\begin{picture}(441.0,108.2)(  4.2,485.1)
  \GGGputW(213.20,537.21){$y$}
  \GGGputW(122.38,579.24){$p_y$}
  \GGGputW[   0, 100](  7.11,589.24){\small non-symplectic method}
  \GGGputW[   0, 100](184.41,589.24){\small $h= 0.05$}
  \GGGputW(438.08,537.21){$y$}
  \GGGputW(347.14,579.24){$p_y$}
  \GGGputW[   0, 100](231.87,589.24){\small St\"ormer--Verlet}
  \GGGputW[   0, 100](409.17,589.24){\small $h= 0.05$}
  \GGGputW(269.33,552.99){\scriptsize zoom$\times 100$}
 \end{picture}
\vspace{-2mm}
\caption{Phase portraits in the $(y,p_y)$-plane for the flat case $\beta=0$. 
\label{fig:figpaper2}}
\vspace{-1mm}
\end{figure}

As noticed in Sect.\ts\ref{sect:martinet},
the normal geodesics in 
the flat case are determined by a one-degree of
freedom Hamiltonian system in the variables $y$ and $p_y$.
We therefore show in Figure~\ref{fig:figpaper2}
the projection onto the $(y,p_y)$-space of
the solutions previously computed
with step size $h=0.05$.
The exact solution starts at $(0,\sin\theta_0 )$ above the
saddle point, turns around the positive stationary point, crosses
the $p_y$-axis at $(0,-\sin\theta_0 )$, turns around the negative
stationary point, and then continues periodically.
The numerical approximation by the
non-symplectic method covers more than one and a half
periods, whereas the \verlet and the exact solution
cover less than one period for the time interval $[0,9]$.
Since the conjugate point is not very sensible with respect to
perturbations in the initial value for $p_y$,
the $(y,p_y)$ coordinates of the
conjugate point obtained by the non-symplectic
integrator are rather accurate, but the corresponding
integration time
is completely wrong.

\begin{table}[b] 
\caption{Accuracy for the first conjugate time.
\label{table:conjtime}}
%\vspace{-1mm}
\centering
\begin{tabular}{|c|c|c|}
\hline
\multicolumn{3}{|c|}{Martinet flat case} \\
\hline
$h$ & \textsc{rk2} & Verlet \\
\hline
$10^{-1}$ & $4.504945$ & $\underline{8}.504716$ \\
$10^{-2}$ & $6.748262$ & $\underline{8.416}622$\\
$10^{-3}$ & $\underline{8.3}60340$ & $\underline{8.4164}12$\\
$10^{-4}$ & $\underline{8.416}349$ & $\underline{8.41641}0$\\
\hline
\multicolumn{3}{|c|}{exact solution: $t_1 \approx \underline{8.416409}$} \\
\hline
\end{tabular}
\qquad\qquad
\begin{tabular}{|c|c|c|}
\hline
\multicolumn{3}{|c|}{Non integrable situation} \\
\hline
$h$ & \textsc{rk2} & Verlet \\
\hline
$10^{-1}$ & $4.511294$ & $\underline{4.88}3832$ \\
$10^{-2}$ & $7.380322$ & $\underline{4.877}056$ \\
$10^{-3}$ & $\underline{4.877}183$ & $\underline{4.87699}8$ \\
$10^{-4}$ & $\underline{4.876997}$ & $\underline{4.876997}$\\
\hline
\multicolumn{3}{|c|}{exact solution: $t_1 \approx \underline{4.876997}$} \\
\hline
\end{tabular}
\end{table}

Table \ref{table:conjtime} lists the
conjugate time obtained with the two integrators using
various step sizes.
There is a significant difference between the two methods. 
We can see that with the \verlet method (\ref{eq:verlet})
a step size of order $h=10^{-2}$ provides
a solution with $4$ correct digits. A step size a
$100$ times smaller is needed 
to get the same precision with the non-symplectic method.

\subsection{Non integrable perturbation} \label{sect:nonintcase}

\begin{figure}[ht]
\centering
%\GGGinput{}{figmartinet2b}
%\GGGinput{../prog/figsb/}{figmartinet2b}
 \begin{picture}(0,0)
  \epsfig{file=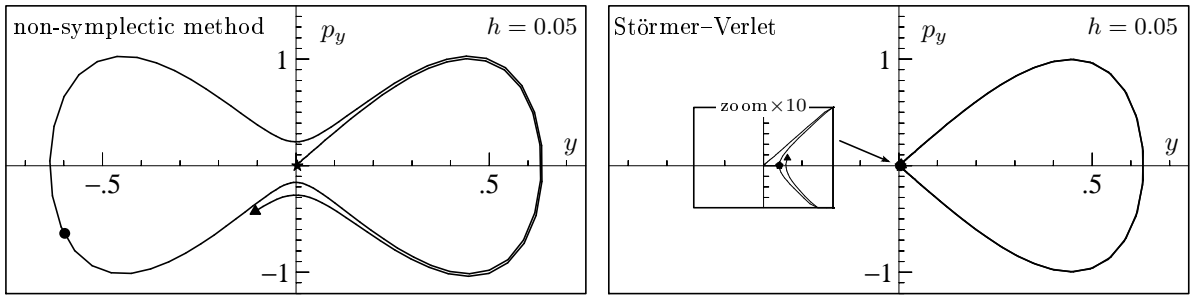}
 \end{picture}%
\begin{picture}(441.0,108.2)(  4.2,485.1)
  \GGGputW(213.20,537.21){$y$}
  \GGGputW(122.38,579.24){$p_y$}
  \GGGputW[   0, 100](  7.11,589.24){\small non-symplectic method}
  \GGGputW[   0, 100](184.41,589.24){\small $h= 0.05$}
  \GGGputW(438.08,537.21){$y$}
  \GGGputW(347.14,579.24){$p_y$}
  \GGGputW[   0, 100](231.87,589.24){\small St\"ormer--Verlet}
  \GGGputW[   0, 100](409.17,589.24){\small $h= 0.05$}
  \GGGputW(271.49,552.99){\scriptsize zoom$\times 10$}
 \end{picture}
\caption{Phase portraits in the $(y,p_y)$-plane for the
non integrable case $\beta=-10^{-4}$.
\label{fig:figpaper2b}}
\end{figure}

For our next numerical experiment we choose the
perturbation parameter $\beta=-10^{-4}$ in the
differential equation (\ref{eq:martinetsyst}).
We consider the same initial values and the same
integration interval as in Sect.\ts\ref{sect:flatcase}.
The exact solution is no longer periodic and,
due to the fact that $\beta$ is chosen negative,
its projection onto the $(y,p_y)$-space
slowly spirals inwards around the positive
stationary point (see right picture in Figure \ref{fig:figpaper2b}).

\begin{figure}[ht]
\centering
%\GGGinput{}{figmartinet1b}
%\GGGinput{../prog/figsb/}{figmartinet1b}
 \begin{picture}(0,0)
  \epsfig{file=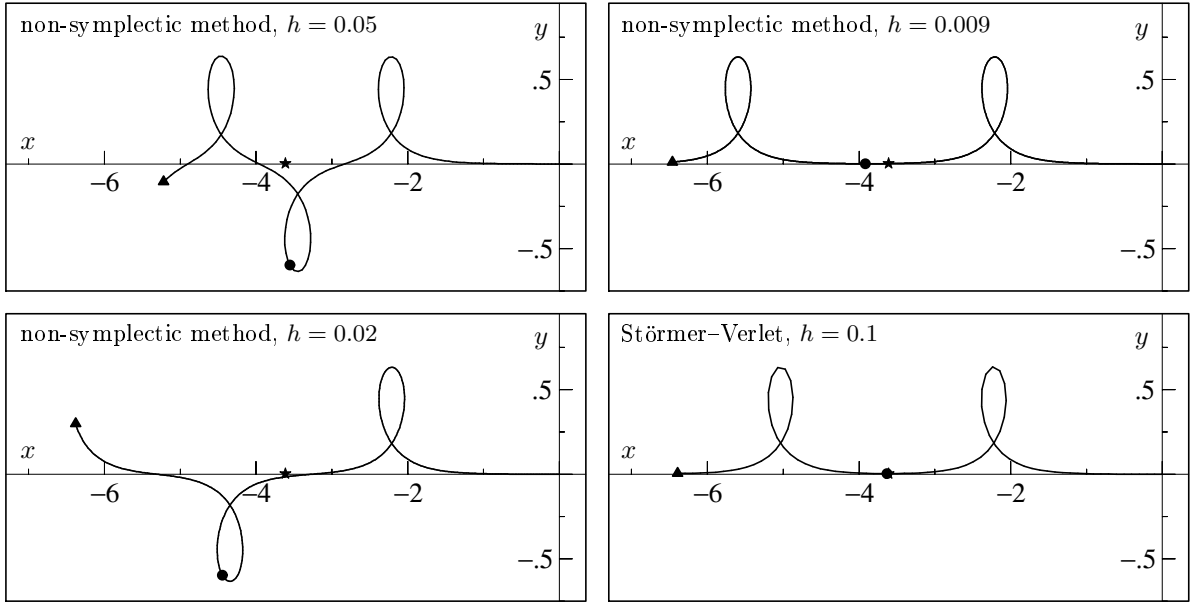}
 \end{picture}%
\begin{picture}(441.0,224.8)(  4.2,368.5)
  \GGGputW[   0, 100](  9.88,543.95){$x$}
  \GGGputW[   0, 100](202.12,586.83){$y$}
  \GGGputW[   0, 100](  9.88,590.08){\small non-symplectic method,  $h= 0.05$}
  \GGGputW[   0, 100](  9.88,427.24){$x$}
  \GGGputW[   0, 100](202.12,470.24){$y$}
  \GGGputW[   0, 100](  9.88,473.37){\small non-symplectic method,  $h= 0.02$}
  \GGGputW[   0, 100](234.64,543.95){$x$}
  \GGGputW[   0, 100](426.87,586.83){$y$}
  \GGGputW[   0, 100](234.64,590.08){\small non-symplectic method,  $h= 0.009$}
  \GGGputW[   0, 100](234.64,427.24){$x$}
  \GGGputW[   0, 100](426.87,470.24){$y$}
  \GGGputW[   0, 100](234.64,473.37){\small St\"ormer--Verlet,  $h= 0.1$}
 \end{picture}
\caption{Trajectories in the $(x,y)$-plane for the non
integrable case $\beta=-10^{-4}$.
\label{fig:figpaper1b}}
\medskip
\end{figure}

Figures \ref{fig:figpaper1b} and \ref{fig:figpaper2b}
and Table \ref{table:conjtime}
display the numerical results obtained by the two
integrators for the differential equation~(\ref{eq:martinetsyst})
with $\beta = -10^{-4}$. 
The interpretation of the symbols (triangles,
circles, and stars) is the same as before.
The excellent behavior of the symplectic integrator is even
more spectacular than 
in the flat case, and the pictures obtained for the
\verlet method agree extremely well
with the exact solution.
The non-symplectic method gives qualitatively wrong
solutions for
step sizes larger than $h=0.01$. In the $(y,p_y)$-space
it alternatively spirals around the right and left
stationary points whereas the exact
solution spirals only around the positive stationary point.
In contrast to the Martinet flat case, the conjugate point
obtained by the non-symplectic method is here wrong also
in the $(y,p_y)$-space.

\subsection{An asymptotic formula on the first conjugate time
in the 
Martinet flat case}\label{sect:conjecture}

Now that we have shown the efficiency of symplectic integrators, we can
make more precise
the asymptotic behaviour studied in \cite{SR97}.
For the initial values of (\ref{eq:init}) and $\beta =0$,
consider the ratio
\begin{equation*}
R=\frac {t_1 \sqrt{p_z}}{3K(k)} ,
\end{equation*}
where $t_1$ is the first conjugate time for the normal geodesic,
and $K(k)$ is an 
elliptic integral of the first kind,
\begin{equation*}
K(k)= \displaystyle \int_0^{\pi/2} \frac1 {\sqrt{1-k^2\sin^2 u}} du, \qquad \quad k=\sin(\theta_0 /2) .
\end{equation*}
By studying analytic solutions for the normal geodesics, it is proved 
in \cite{SR97} that this ratio
satisfies the inequality $2/3 \leq R \leq 1$.
It follows from a rescaling of the equations (\ref{eq:martinetsyst})
that $R$ is independent of $p_z$.

\begin{figure}[ht]
\centering
%\GGGinput{}{figconj}
%\GGGinput{../prog/figconj/}{figconj}
 \begin{picture}(0,0)
  \epsfig{file=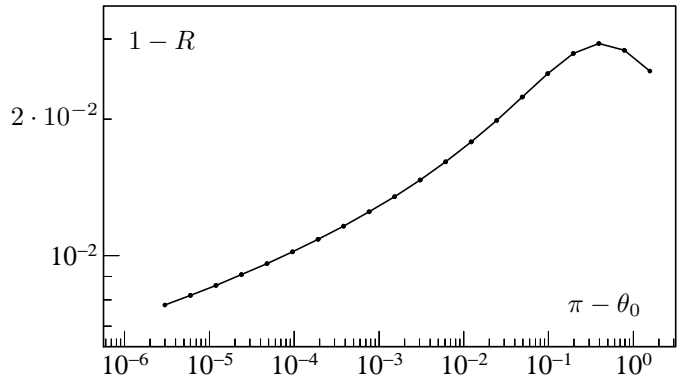}
 \end{picture}%
\begin{picture}(213.3,128.0)(  0.0,469.4)
  \GGGputW[   0, 100](-34.69,561.66){$2 \cdot 10^{-2}$}
  \GGGputW[   0, 100](  9.40,589.00){$1-R$}
  \GGGputW[   0, 100](173.81,489.87){$\pi-\theta_0$}
 \end{picture}
\caption{Illustration of the asymptotic behaviour of $R$
(\verlet scheme with step size $h=10^{-4}$).
\label{fig:conj}}
\end{figure}

In Figure \ref{fig:conj}, we represent the values of $1-R$ as a function of 
$\varepsilon = \pi-\theta_0$, for various initial values
$\theta_0$. 
The numerical results indicate that the ratio $R$ depends on $\theta_0$, and 
$R \longrightarrow 1^-$ slowly for $\theta_0 \longrightarrow \pi^-$.

\section{Backward error analysis}
\label{sect:backward}

The theory of backward error analysis is fundamental for the study 
of geometric integrators and it is treated in much detail in the 
monographs of Sanz-Serna \& Calvo \cite{sanz-serna94nhp},
Hairer, Lubich \& Wanner \cite[Chap.\ts IX]{hairer06gni}, 
and Leimkuhler \& Reich \cite{leimkuhler04shd}.
It allows us to explain the numerical phenomena
encountered in the previous section.

\subsection{Backward error analysis and energy conservation}\label{sect:bea1}

We briefly present the main ideas of backward error analysis for the 
study of symplectic integrators, see \cite[Chap.\ts IX]{hairer06gni}.
Consider a system of differential equations
\begin{equation} \label{eq:eqf}
\dot y = f(y), \quad y(0)=y_0
\end{equation}
and a numerical integrator $y_{n+1}=\Phi_h(y_n)$ of order $p$. 
The idea is to search for a \textit{modified differential equation} 
written as a formal series in powers of the step size $h$,
\begin{equation} \label{eq:modifeq}
\dot {\widetilde y} = 
\widetilde f(\widetilde y) = 
f(\widetilde y) + h^pf_{p+1}(\widetilde y) + h^{p+1}f_{p+2}(\widetilde y) + \ldots ,
\end{equation}
such that $y_n=\widetilde y(t_n)$ for $t_n=nh, \, n=0,1,2,\ldots$, in the 
sense of formal power series. The motivation of this approach is that it is often easier to study the 
modified equation (\ref{eq:modifeq}) than
directly the numerical solution. 

What makes backward error analysis so important for the study of symplectic 
integrators is the fact that, when applied to a Hamiltonian system 
$\dot y = J^{-1} \nabla H(y)$, the modified equation (\ref{eq:modifeq}) 
has the same structure
$\dot {\widetilde y} = J^{-1} \nabla \widetilde H(\widetilde y)$
with a \textit{modified Hamiltonian} 
\begin{equation*}
\widetilde H(y) = H(y) + h^pH_{p+1}(y) +  h^{p+1}H_{p+2}(y) +\ldots ~.
\end{equation*}
However, the series usually diverges, so a truncation at a suitable 
order~$N(h)$ is necessary, 
\begin{equation*}
\widetilde H(y) = H(y) + h^pH_{p+1}(y) + \ldots + h^{N-1}H_{N}(y)  .
\end{equation*}
This truncation induces an error that can be made exponentially small,
by choosing $N(h)\sim C/h$,
see \cite[Theorem IX.8.1]{hairer06gni}.
More precisely, we have that for $t_n=nh$ and $h \rightarrow 0$,
\begin{equation} \label{eq:consHtilde} 
\widetilde H(y_n) = \widetilde H(y_0) + \bigo(t_n e^{- h_0 / h}) .
\end{equation}
as long as the numerical solution $\{y_n\}$ stays in a compact set.
On intervals of length
$\bigo(e^{h_0 / 2h})$, the modified Hamiltonian
$\widetilde H(y )$ is thus exactly conserved up to
exponentially small terms.

\subsection{Backward error analysis for the Martinet problem}

Symplectic integrators are successfully applied in the
long-time integration 
of Hamiltonian systems, for instance
in astronomy (e.g. the Outer Solar System over 100 million years 
\cite[Sect.\ts I.2.4]{hairer06gni}), or in molecular dynamics \cite[Chap.\ts 11]{leimkuhler04shd}.
Here the situation is quite different because we are interested 
in the numerical integration of Hamiltonian systems
on relatively short time intervals.

\subsubsection{Martinet flat case}

Consider the Martinet problem (\ref{eq:martinetsyst}) 
in the flat case $\beta=0$.
Its interesting dynamics takes place in the
$(y,p_y)$ plane, and it is not influenced by the other
variables (only by their initial values).
We put $\eta = (y,p_y)$, and we denote by $f(\eta )$ the
Hamiltonian vector field composed by the
corresponding two equations
of (\ref{eq:martinetsyst}).
For a numerical integrator of order $p=2$, the associated
modified 
differential equation has the form
\begin{equation}\label{eq:modifeq2}
\dot {\widetilde \eta} = f(\widetilde \eta) +
 h^2 f_3 (\widetilde \eta) + \bigo(h^3) .
\end{equation}

Consider first the symplectic \verlet method.
It follows from
Sect.\ts \ref{sect:bea1} that its modified
differential equation is Hamiltonian,
and from (\ref{eq:consHtilde})
that the modified Hamiltonian $\widetilde H(\eta )$ is preserved up
to exponentially small terms along the numerical solution. This implies
that the numerical solution remains exponentially close to a periodic
orbit in the $(y,p_y)$-space. The critical point
$(y\! = \! 0,\, p_y\! = \! 0)$ is a saddle point also for the
modified differential equation
(because the origin is stationary also for the numerical
solution and thus for the modified equation).
Therefore, any numerical solution starting close to the
origin has to come back to it after turning around one of
the stationary points. The minimal distance to the origin
will always stay the same (see the zoom in
Figure \ref{fig:figpaper2}). This explains the good behavior
of symplectic integrators.

For the non-symplectic integrator, the term $h^2f_3(\eta )$ is not
Hamiltonian. Therefore the solution of the modified differential equation
(and hence also the numerical solution) is no longer periodic. In fact,
it spirals outwards and after surrounding the first
stationary point, the numerical solution does not approach
the saddle point sufficiently close, which induces a faster dynamics
as can be observed in Figures \ref{fig:figpaper1} and
\ref{fig:figpaper2}.
This causes a huge error, because close to the saddle point
the numerical solution is most sensible to errors.

\subsubsection{Non integrable perturbation}
\label{sect:backnonint}

In this case, the argument in the comparison of
symplectic and non-symplectic integrators 
is very similar to the discussion of the Van der Pol's equation
in \cite[Sect.\ts XII.1]{hairer06gni}. 
For $\beta\neq 0$ (non integrable perturbation), the dynamics
takes place in the four dimensional space with variables
$\eta =(x,y,p_x,p_y)$. In this space
the system (\ref{eq:martinetsyst}) becomes
\begin{equation*}
\dot \eta = f(\eta ) + \beta g(\eta )
\end{equation*}
where $f(\eta )$ is the Hamiltonian vector field corresponding
to $\beta =0$ and
$g(y)=\bigo(1)$ depends smoothly on $\beta$. 
Here, the modified equation becomes
\begin{equation*}
\dot {\widetilde \eta} = f(\widetilde \eta ) +\beta g(\widetilde \eta )
 + h^2 f_3 (\widetilde \eta ) + \bigo(h^3 +\beta h^2) ,
\end{equation*}
where the perturbation term $h^2 f_3 (\eta )$
is the same as for the Martinet flat case.

For the symplectic integrator,
the perturbation $\beta g(\eta )$ has the same effect
for the original problem as for
$\dot {\widetilde \eta} = f(\widetilde \eta )
 + h^2 f_3 (\widetilde \eta ) +\ldots \, $. This
explains the correct qualitative behavior 
for small $h$ and small $\beta$. There is no restriction on the
step size $h$ compared to the size of $\beta$. 

For the non-symplectic integrator, each of the perturbation terms
$\beta g(\eta )$ and $h^2 f_3 (\eta ) $
destroys the periodic orbits in the subsystem for the
$(y,p_y)$ variables,
and the dominant one will determine the behavior of the numerical
solution. Only when $h^2 \ll |\beta |$, the numerical solution will
catch the correct dynamics of the problem.
In Figures \ref{fig:figpaper1b} and \ref{fig:figpaper2b}, where $\beta =-10^{-4}$,
this condition is not satisfied for $h\ge 10^{-2}$.
Since $\beta$ is chosen small and
negative, the two perturbation terms are conflicting.
The term
$\beta g(\eta )$ causes the solution to spiral around the
positive stationary point,
whereas the term $h^2 f_3 (\eta ) $ causes it to spiral
alternatively around both stationary points.
For too large step sizes the qualitative behavior of the non-symplectic
integrator (\ref{eq:rk2}) is thus completely wrong.

\begin{rmrk} \rm
The problem (\ref{eq:martinetsyst})
with $\beta =0$ has a lot of
symmetries. In the $(y,p_y)$-space the orbits are
symmetric with respect to the $y$-axis and also
with respect to the $p_y$-axis. If we apply a symmetric
numerical integrator (not necessarily symplectic), it is
possible to prove the same qualitative behavior as for
the symplectic \verlet method.
This follows from the fact that the solution
of the modified equation (numerical orbit) corresponding to
a symmetric method has the same symmetry properties
as the exact flow
(see \cite[Sect.\ts IX.2]{hairer06gni} for precise statements).
Consequently, in the $(y,p_y)$ plane and for $\beta =0$,
the solution will stay exponentially near to a closed orbit,
as it is the case for symplectic integrators. In the non
integrable case, the good behavior of symmetric methods
can be explained as in Sect.\ts
\ref{sect:backnonint} for symplectic methods.
\end{rmrk}

%%-----------------------------
%%      your bibliography
\bibliographystyle{plain}

%\newpage
\tableofcontents

\end{document}